\documentclass[siam10.clo]{siamltex}

\newcommand{\N}{I\!\! N}

\begin{document}
\vspace*{-40pt}
\rm

\vspace*{16pt}

\vspace{1.2cm}

\def\bfE{\mbox{\boldmath$E$}}
\def\bfG{\mbox{\boldmath$G$}}

\title{Quadratic categories, Koszul resolutions.}

\author{Aristide Tsemo, \thanks{ Department of Mathematics and Computer Science,
        350, Victoria Street
        Toronto, ON
        M5B 2K3
        Canada  ({\tt taristid@scs.ryerson.ca}).}
        \and  Isaac Woungang\thanks{Department of Mathematics
        and Computer Science,
        350 Victoria Street
        Toronto, ON
        M5B 2K3
        Canada
        {\tt iwoungan@scs.ryerson.ca }}}

\maketitle

\begin{abstract}

The category of quadratic algebras has been endowed by Manin [6]
with two tensor products. These products have been generalized to
quadratic operads by Ginzburg and Kapranov [5], and to
$n$-homogeneous algebras by Berger [2]. The purpose of this paper is
to define an abstract notion of quadratic category such that the
categories of quadratic algebras and quadratic operads are examples
of this notion. We define Koszul complexes in this setting,
representations of quadratic categories in the category of quadratic
algebras, and Tannakian quadratic categories.
 \end{abstract}

\section{Introduction}

 A quadratic algebra is a quotient of a tensor algebra
$T(V)$ of a finite dimensional vector space $V$, by an ideal $C$
generated by a subspace of $V\otimes V$. These algebras are
${\N}$-graded algebras generated by their elements of degree $1$,
which satisfy quadratic relations. Quadratic algebras appear in
different domains of mathematics, as in topology with the notion of
Steenrood algebras, in differential geometry, the Clifford algebras
are one of the main tools to study $Spin$-geometry, in group theory,
symmetric and exterior algebras are very useful. In algebraic
geometry, a projective scheme can be realized as  a projective
spectrum of a quadratic algebra, this is equivalent to saying that
the category of quasi-coherent sheaves over a projective scheme $N$,
is equivalent to the category of graded modules  over a quadratic
algebra.

Cohomology theories are defined in the general context of abelian
categories with enough injective objects by applying the $Hom$
functor to resolutions. To compute cohomology groups, we need to
define complexes which represent these   resolutions, like the
Chevalley complex in   groups theory, the Koszul complex in Lie
algebras theory, the Bar complex in  associative algebras theory...
Even at this stage, these canonical complexes are hardly tractable.
This has motivated Priddy in [7] to define a generalized Koszul
complex for quadratic algebras which allows to compute their
cohomology when it is a resolution. The complex defined by Priddy is
a generalization of the classical Koszul complex defined with the
exterior and symmetric algebras. This is useful in practice since
the Koszul complex is simpler than the Bar resolution.

 The automorphisms group of quadratic algebras has been used in
 theoretical physics in the inverse scattering problem, and in low
 dimensional topology. It is in this context that
  Manin has endowed in [6] the category of quadratic algebras
with two tensor products $o$ and $\bullet$, he has also defined
the notion of quadratic dual which represents the Yoneda algebra
of quadratic Koszul algebras. He has shown that there exists
internal $\underline{Hom}$ object in the category of quadratic
algebras endowed with the tensor product $\bullet$.

An operad is an object which encodes operations. These objects have
been defined in homotopy theories, and  in [5] to study algebraic
structures. The Manin tensor products have been adapted to the
theory of quadratic operads by Ginzburg and Kapranov [5], in their
paper they have defined the notion of Koszul resolution of a
quadratic operad, which is an application of the general Koszul
duality defined by Beilinson Ginzburg and Schechtman [1]. Recently,
Berger [3] has defined the category of $n$-Koszul algebras. He has
endowed this category with  tensor products similar to the ones
defined by Manin.

The purpose of this paper is to define a general notion of quadratic
category endowed with two tensor products and a duality which
satisfy some compatibility conditions.  The categories of quadratic
algebras, quadratic operads, and $n$-Koszul algebras are examples of
quadratic categories. In this context, we show the following result
as theorem 2.2:

\bigskip
{\it

 Let $(C,o,\bullet,!)$ be a quadratic category, then for each
objects $U$ and $V$ of $C$, $V o U^!$ is an internal $Hom$ of the
tensor category $(C,\bullet)$.}
\bigskip

Tensor categories have been studied in  [4] to determine properties
of the cohomology ring of algebraic varieties. These authors have
defined the notion of Tannakian category which is an Abelian rigid
tensor category endowed with an exact faithful functor to a category
of vector spaces, and have shown that a Tannakian category is
equivalent to the category of representations of an affine group
scheme. We adapt the study of these authors to  quadratic categories
by defining a quadratic Tannakian category: it is  an Abelian
quadratic category endowed with an exact faithful functor to the
category of quadratic algebras, we show: (see theorem 4.5)
\bigskip

{\it  A quadratic Tannakian category is equivalent to the category
of quadratic representations of an affine group scheme. }

\medskip

A projective scheme $U$, is the projective spectrum of a quadratic
algebra $T_U$, we can associate to $U$ the quadratic Tannakian
category $T'_U$ generated by $T_U$. The property of the group $H_U$
whose category of quadratic representations is equivalent to $T'_U$
seems to by an interesting object to study.  Many authors have tried
to define non commutative algebraic geometry. The algebraic geometry
of tensor categories is studied by Deligne,  perhaps quadratic
Tannakian categories represent the  good framework for
noncommutative algebraic geometry, and the category of quadratic
algebras the motivic category in non commutative algebraic geometry.

\section{ Quadratic categories}

The purpose of this part is to present the notion and properties
of quadratic categories.

\begin{definition}

A quadratic category  $(C,o,\bullet)$ is a category $C$, endowed
with two tensors products $o$ and $\bullet$, whose neutral elements
are respectively $I_o$ and $I_{\bullet}$, and whose associativity
constraints are respectively $c_o$ and $c_{\bullet}$. We suppose
that the following properties are satisfied: For each objects $U,
U_1, U_2, U_3$ of $C$, there exists an object $U^!$ of $C$,
morphisms

$c_U:I_{\bullet}\rightarrow U o U^!$, $d_U:U^!\bullet U\rightarrow
I_0$, $f_{U_1,U_2,U_3}: (U_1 o U_2)\bullet U_3\longrightarrow U_1 o
(U_2 \bullet U_3)$, $h_{U_1,U_2,U_3}:U_1\bullet (U_2 o
U_3)\longrightarrow (U_1 \bullet U_2) o U_3$

such that the following diagrams are commutative:

\begin{equation}
\matrix{(U_1\bullet (U_2 o U_3))\bullet U_4 &{\buildrel
{c_{\bullet}(U_1,U_2 o U_3, U_4)}\over{\longrightarrow}} U_1\bullet
((U_2 o U_3)\bullet U_4){\buildrel {Id_{U_1}\bullet
f_{U_2,U_3,U_4}}\over {\longrightarrow}} & U_1\bullet (U_2 o
(U_3\bullet U_4))\cr \downarrow h_{U_1, U_2, U_3}\bullet Id_{U_4}&&
\downarrow h_{U_1,U_2, U_3\bullet U_4} \cr ((U_1\bullet U_2) o
U_3)\bullet U_4  &{\buildrel{f_{U_1\bullet
U_2,U_3,U_4}}\over{\longrightarrow}}& (U_1\bullet U_2) o (U_3\bullet
U_4)}
\end{equation}

\begin{equation}
\matrix{(U_1 o U_2)\bullet (U_3 o U_4)&{\buildrel{f_{U_1,U_2,U_3 o
U_4}}\over{\longrightarrow}} &U_1 o (U_2\bullet(U_3 o U_4))\cr
\downarrow h_{U_1 o U_2,U_3,U_4}&&\downarrow Id_{U_1} o
h_{U_2,U_3,U_4}\cr ((U_1 o U_2)\bullet U_3) o U_4
&{\buildrel{f_{U_1,U_2,U_3}o Id_{U_4}}\over{\longrightarrow}} (U_1
o(U_2\bullet U_3)) o U_4 {\buildrel{c_0(U_1,U_2\bullet
U_3,U_4)}\over{\longrightarrow}} & U_1 o ((U_2\bullet U_3) o U_4)}
\end{equation}

\begin{equation}
\matrix{I_{\bullet}\bullet U&{\buildrel{c_U\bullet
Id_U}\over{\longrightarrow}} (U o U^!)\bullet
U{\buildrel{f_{U,U^!,U}}\over{\longrightarrow}} U o (U^!\bullet
U){\buildrel{Id_U o d_U}\over{\longrightarrow}}& U o I_o\cr
\downarrow &&\downarrow \cr
U&{\buildrel{Id_U}\over{\longrightarrow}} &U}
\end{equation}

\begin{equation}
\matrix{U^!\bullet I_{\bullet}&{\buildrel{Id_{U^!}\bullet
c_U}\over{\longrightarrow}} U^! \bullet (U o
U^!){\buildrel{h_{U^!,U,U^!}}\over{\longrightarrow}} (U^!\bullet U)
o U^!{\buildrel{{d_U} o Id_{U^!}}\over{\longrightarrow}}& I_o o
U^!\cr \downarrow &&\downarrow \cr
U^!&{\buildrel{Id_{U^!}}\over{\longrightarrow}} &U^!}
\end{equation}

Let $u_1:U_1\rightarrow U'_1$, $u_2:U_2\rightarrow U'_2$ and
$u_3:U_3\rightarrow U'_3$ be three morphisms of $C$. The following
diagrams are supposed to be commutative:

\begin{equation}
\matrix{U_1\bullet (U_2 o U_3)&{\buildrel{u_1\bullet (u_2 o
u_3)}\over{\longrightarrow}}& U'_1\bullet (U'_2 o U'_3)\cr
\downarrow h_{U_1,U_2,U_3}&&\downarrow h_{U'_1,U'_2,U'_3}\cr
(U_1\bullet U_2) o U_3 &{\buildrel{(u_1\bullet u_2)o
u_3}\over{\longrightarrow}}&(U'_1\bullet U'_2) o U'_3}
\end{equation}

\begin{equation}
\matrix{(U_1 o U_2)\bullet U_3 &{\buildrel{(u_1 o u_2)\bullet
u_3}\over{\longrightarrow}} &(U'_1 o U'_2)\bullet U'_3\cr\downarrow
f_{U_1,U_2,U_3}&&\downarrow f_{U'_1,U'_2,U'_3}\cr U_1 o (U_2\bullet
U_3)&{\buildrel{u_1 o (u_2\bullet u_3)}\over{\longrightarrow}}& U'_1
o (U'_2\bullet U'_3)}
\end{equation}
\end{definition}

The example who has motivated the construction of quadratic
categories is the theory of quadratic algebras, recall the
constructions of the tensor products defined by Manin in [6].

Let $U=T(U_1)/(C_1)$ and $V=T(V_1)/(C_2)$ be two quadratic algebras
respectively isomorphic to the quotient of the tensor algebra of the
finite dimensional $L$-vector spaces $U_1$ and $V_1$, by the ideals
generated by $C_1\subset U_1\otimes U_1$ and $C_2\subset V_1\otimes
V_1$. We endow the class of quadratic algebras with the structure of
a category such that $Hom(U,V)$ is the set of morphisms of algebras
$h:U\rightarrow V$ defined by a linear application
$h_1:U_1\rightarrow V_1$ such that $(h_1\otimes h_1)(C_1)\subset
C_2$.

The quadratic algebra $U\bullet V$ is the quotient of the tensor
algebra  $T(U_1\otimes V_1)$ by the ideal generated by
$t_{23}(C_1\otimes C_2)$. The isomorphism $t_{23}:$
$U_1^{\otimes^2}\otimes V_1^{\otimes^2}\rightarrow(U_1\otimes
V_1)^{\otimes^2}$ is defined by $t_{23}(u_1\otimes u'_1\otimes
v_1\otimes v'_1)=u_1\otimes v_1\otimes u'_1\otimes v'_1$.

The tensor product $U o V$ is  the quotient of the tensor algebra
$T(U_1\otimes V_1)$ by the ideal generated by
$t_{23}(U_1^{\otimes^2}\otimes C_2+  C_1\otimes V_1^{\otimes^2})$.

The quadratic dual $U^!$ of $U$ is the quadratic algebra
$T(U_1^*)/C_1^*$, where $U_1^*$ is the dual vector space of $U_1$,
and $C_1^*$ the annihilator of $C_1$ in $U^*_1\otimes U^*_1$.

The neutral element $I_o$ of $o$ is $T(L)$, and the neutral element
of $\bullet$,  $I_{\bullet}$ is $L$. The map
$c_U:I_{\bullet}\rightarrow U o U^!$ is defined by the map
$c'_U:L\rightarrow U_1\otimes U_1^*$, $c'_U(1)\rightarrow
\sum_{i=1}^{i=n}u_i\otimes u^i$, where $(u_1,..,u_n)$ is a basis of
$U_1$,  and $(u^1,..,u^n)$ its dual basis. The map $d_U:U^!\bullet
U\rightarrow I_o$ is defined by the duality $U_1\otimes
U_1^*\rightarrow L$.

Let $U^i=T(U_i)/(C_i), i=1,2,3$. The ideal $C$ which defines $(U^1 o
U^2)\bullet U^3$ is $t_{23}(t_{23}(U_1^{\otimes^2}\otimes C_2 +
C_1\otimes U_2^{\otimes^2})\otimes C_3)$, and the ideal $C'$ which
defines $U^1 o (U^2\bullet U^3)$ is
$t_{23}(U_1^{\otimes^2}\otimes(t_{23}(C_2\otimes C_3))+C_1\otimes
(U_2\otimes U_3)^{\otimes^2}$). We remark that $C$ is contained in
$C'$, the associativity constraint $(U_1\otimes U_2)\otimes
U_3\longrightarrow U_1\otimes (U_2\otimes U_3)$ of vector spaces
projects to define the quadratic constraint $f_{U^1,U^2,U^3}$.

The ideal $D$ which defines $U^1\bullet (U^2 o U^3)$ is
$t_{23}(C_1\otimes t_{23}(U_2^{\otimes^2}\otimes C_3 + C_2\otimes
U_3^{\otimes^2}))$. The ideal $D'$ which defines the algebra
$(U^1\bullet U^2) o U^3$ is $t_{23}({(U_1\otimes
U_2)}^{\otimes^2}\otimes C_3+t_{23}(C_1\otimes C_2)\otimes
U_3^{\otimes^2})$. We remark that $D'$ contains $D$ this implies
that the associativity constraint for vector spaces $U_1\otimes
(U_2\otimes U_3)\longrightarrow (U_1\otimes U_2)\otimes U_3$
projects to define the quadratic associativity constraint
$h_{U^1,U^2,U^3}$.

\begin{theorem}

Let $C$ be a quadratic category,  $L$, and $N$ two objects of $C$.
The contravariant functor:

\begin{eqnarray*}
C\longrightarrow C
\end{eqnarray*}

\begin{eqnarray*}
U\longrightarrow Hom(U\bullet L,N)
\end{eqnarray*}

is representable by $N o {L}^!$.
\end{theorem}

\begin{proof}

We define a map between the functors $U\rightarrow Hom(U\bullet
L,N)$ and $U\rightarrow Hom(U,N o L^!)$ by assigning to each element
$u$ in $Hom(U\bullet L,N)$, the element $u'$ in $Hom(U,N o L^!)$
defined by:

\begin{eqnarray*}
\matrix{U\rightarrow U\bullet I_{\bullet}{\buildrel{Id_U\bullet
c_L}\over{\rightarrow}} U\bullet(L o
L^!){\buildrel{h_{U,L,L^!}}\over {\longrightarrow}} (U\bullet L) o
L^!{\buildrel{u o Id_{L^!}}\over{\longrightarrow}} N o L^!}
\end{eqnarray*}

We define a map between the functors $U\rightarrow Hom(U,N o L^!)$
and $U\rightarrow Hom(U\bullet N,L)$ by assigning to each element
$v$ in $Hom(U,N o L^!)$ the element $v"$ in $Hom(U\bullet L,N)$
defined by:

\begin{eqnarray*}
\matrix{U\bullet L {\buildrel{v\bullet
Id_L}\over{\longrightarrow}}(N o L^!)\bullet
L{\buildrel{f_{N,L^!,L}}\over{\longrightarrow}} N o (L^!\bullet
L){\buildrel{Id_N o d_L}\over{\longrightarrow}}}N
\end{eqnarray*}

The correspondence $u\rightarrow u'$ and $v\rightarrow v"$ are
morphisms of functors since they are compositions of morphisms of
functors. We have to show that the correspondence defined on
$Hom(U\bullet L,N)$ by $u\rightarrow (u')"$ is  the identity
morphism of the functor $U\rightarrow Hom(U\bullet L,N)$. We have:

\begin{eqnarray*}
(u')"=\matrix{(U\rightarrow U\bullet
I_{\bullet}{\buildrel{Id_U\bullet c_L}\over{\rightarrow}} U\bullet(L
o L^!){\buildrel{h_{U,L,L^!}}\over {\longrightarrow}} (U\bullet L) o
L^!{\buildrel{u o Id_{L^!}}\over{\longrightarrow}} N o L^!)\bullet
L}
 \matrix{{\buildrel{f_{N,L^!,L}}\over{\longrightarrow}} N o
(L^!\bullet L){\buildrel{Id_N o d_L}\over{\longrightarrow}}}N
\end{eqnarray*}

Let $(C,\otimes)$ be a tensor category, and $u:U\rightarrow U'$,
$u':U'\rightarrow U"$, $v:V\rightarrow V'$, $v":V'\rightarrow V"$
arrows of $C$, $(u'\otimes v')\circ (u\otimes v)=(u'\circ u)\otimes
(v'\circ v)$. Applying this fact to the tensor product $\bullet$ at
the first line of the previous equality, we obtain:

\begin{eqnarray*}
(u')"=\matrix{U\bullet L\rightarrow (U\bullet I_{\bullet})\bullet
L{\buildrel{(Id_U\bullet c_L)\bullet Id_L}\over{\rightarrow}}
(U\bullet(L o L^!))\bullet L{\buildrel{h_{U,L,L^!}\bullet Id_L}\over
{\longrightarrow}} ((U\bullet L) o L^!)\bullet L{\buildrel{(u o
Id_{L^!})\bullet Id_L}\over{\longrightarrow}} (N o L^!)\bullet L}
\end{eqnarray*}

\begin{eqnarray*}
\matrix{{\buildrel{f_{N,L^!,L}}\over{\longrightarrow}} N o
(L^!\bullet L){\buildrel{Id_N o d_L}\over{\longrightarrow}}N}
\end{eqnarray*}

 Applying property $(2.6)$ we obtain
 \begin{eqnarray*}
 \matrix{((U\bullet L) o L^!)\bullet
L{\buildrel{(u o Id_{L^!})\bullet Id_L}\over{\longrightarrow}} (N o
L^!)\bullet L {\buildrel{f_{N,L^!,L}}\over{\longrightarrow}} N o
(L^!\bullet L){\buildrel{Id_N o d_L}\over{\longrightarrow}}N}
\end{eqnarray*}

\begin{eqnarray*}
 =\matrix{((U\bullet L) o L^!)\bullet
L \matrix{{\buildrel{f_{U\bullet L,  L^!,L}}\over
{\longrightarrow}} (U\bullet L) o (L^!\bullet L){\buildrel {u o
Id_{L^!\bullet L}}\over {\longrightarrow}}
 N o (L^!\bullet L){\buildrel{Id_N o d_L }\over{\longrightarrow}}N}}
\end{eqnarray*}

We deduce that
\begin{eqnarray*}
 (u')"=\matrix{U\bullet L\rightarrow
(U\bullet I_{\bullet})\bullet L{\buildrel{(Id_U\bullet c_L)\bullet
Id_L}\over{\rightarrow}} (U\bullet(L o L^!))\bullet
L{\buildrel{h_{U,L,L^!}\bullet Id_L}\over {\longrightarrow}}
((U\bullet L) o L^!)\bullet L}
\end{eqnarray*}

\begin{eqnarray*}
 \matrix{{\buildrel{f_{U\bullet L,
L^!,L}}\over {\longrightarrow}} (U\bullet L) o (L^!\bullet
L){\buildrel {u o Id_{L^!\bullet L}}\over {\longrightarrow}}
 N o (L^!\bullet L){\buildrel{Id_N o d_L }\over{\longrightarrow}}N}
\end{eqnarray*}
 Applying $(2.1)$ we obtain that

\begin{eqnarray*}
\matrix{(U\bullet(L o L^!))\bullet L{\buildrel{h_{U,L,L^!}\bullet
Id_L}\over {\longrightarrow}} ((U\bullet L) o L^!)\bullet L
{\buildrel{f_{U\bullet L,  L^!,L}}\over {\longrightarrow}} (U\bullet
L) o (L^!\bullet L)}
\end{eqnarray*}

\begin{eqnarray*}
 =\matrix{(U\bullet(L o L^!))\bullet
L{\buildrel{{c_{\bullet}}(U,L o L^!,L)}\over{\longrightarrow}}
U\bullet ((L o L^!)\bullet L){\buildrel{Id_U \bullet
f_{L,L^!,L}}\over {\longrightarrow}}U\bullet(L o (L^!\bullet
L)){\buildrel{h_{U,L,L^!\bullet L}}\over {\longrightarrow}}
(U\bullet L) o (L^!\bullet L)}
\end{eqnarray*}

We deduce that

\begin{eqnarray*}
(u')"=\matrix{U\bullet L\rightarrow (U\bullet I_{\bullet})\bullet
L{\buildrel{(Id_U\bullet c_L)\bullet Id_L}\over{\rightarrow}}
(U\bullet(L o L^!))\bullet L{\buildrel{{c_{\bullet}}(U,L o
L^!,L)}\over{\longrightarrow}} U\bullet ((L o L^!)\bullet L)}
\end{eqnarray*}

\begin{eqnarray*}
\matrix{{\buildrel{Id_U \bullet f_{L,L^!,L}}\over
{\longrightarrow}}U\bullet(L o (L^!\bullet
L)){\buildrel{h_{U,L,L^!\bullet L}}\over {\longrightarrow}}
(U\bullet L) o (L^!\bullet L){\buildrel{u o Id_{L^!\bullet L}}\over
{\longrightarrow}} N o (L^!\bullet L){\buildrel {Id_N o
d_L}\over{\longrightarrow}}N }
\end{eqnarray*}

Using the fact that $o$ is a tensor product we deduce that:

\begin{eqnarray*}
(u')"=\matrix{U\bullet L\rightarrow (U\bullet I_{\bullet})\bullet
L{\buildrel{(Id_U\bullet c_L)\bullet Id_L}\over{\rightarrow}}
(U\bullet(L o L^!))\bullet L{\buildrel{{c_{\bullet}}(U,L o
L^!,L)}\over{\longrightarrow}} U\bullet ((L o L^!)\bullet L)}
\end{eqnarray*}

\begin{eqnarray*}
\matrix{{\buildrel{Id_U \bullet f_{L,L^!,L}}\over
{\longrightarrow}}U\bullet(L o (L^!\bullet
L)){\buildrel{h_{U,L,L^!\bullet L}}\over {\longrightarrow}}
(U\bullet L) o (L^!\bullet L){\buildrel{Id_{U\bullet L} o d_L}\over
{\longrightarrow}} (U\bullet L) o I_o {\buildrel
u\over{\longrightarrow}}N }
\end{eqnarray*}
Using property $(2.5)$, we obtain that:

\begin{eqnarray*}
\matrix{U\bullet(L o (L^!\bullet L)){\buildrel{h_{U,L,L^!\bullet
L}}\over {\longrightarrow}} (U\bullet L) o (L^!\bullet
L){\buildrel{Id_{U\bullet L} o d_L}\over {\longrightarrow}}
(U\bullet L) o I_o}
\end{eqnarray*}

\begin{eqnarray*}
 =\matrix{U\bullet (L o (L^!\bullet
L)){\buildrel{Id_U\bullet(Id_L o d_L)}\over{\longrightarrow}}
U\bullet (L o I_o){\buildrel
{h_{U,L,I_o}}\over{\longrightarrow}}(U\bullet L) o I_o}
\end{eqnarray*}

We deduce that

\begin{eqnarray*}
(u')"=\matrix{U\bullet L\rightarrow (U\bullet I_{\bullet})\bullet
L{\buildrel{(Id_U\bullet c_L)\bullet Id_L}\over{\rightarrow}}
(U\bullet(L o L^!))\bullet L{\buildrel{{c_{\bullet}}(U,L o
L^!,L)}\over{\longrightarrow}} U\bullet ((L o L^!)\bullet L)}
\end{eqnarray*}

\begin{eqnarray*}
\matrix{{\buildrel{Id_U \bullet f_{L,L^!,L}}\over
{\longrightarrow}}U\bullet (L o (L^!\bullet
L)){\buildrel{Id_U\bullet(Id_L o d_L)}\over{\longrightarrow}}
U\bullet (L o I_o){\buildrel
{h_{U,L,I_o}}\over{\longrightarrow}}(U\bullet L) o I_o{\buildrel
u\over{\longrightarrow}}N }
\end{eqnarray*}

Using the fact that the associative constraint $c_{\bullet}$ of the
tensor product $\bullet$ is a functorial correspondence, we obtain
that

\begin{eqnarray*}
\matrix{ (U\bullet I_{\bullet})\bullet L{\buildrel{(Id_U\bullet
c_L)\bullet Id_L}\over{\longrightarrow}} (U\bullet(L o L^!))\bullet
L{\buildrel{c_{\bullet}(U,L o L^!,L)}\over{\longrightarrow}}
U\bullet ((L o L^!)\bullet L)}
\end{eqnarray*}

\begin{eqnarray*}
=\matrix{(U\bullet I_{\bullet})\bullet
L{\buildrel{{c_{\bullet}}(U,I_{\bullet},L})\over{\longrightarrow}}U\bullet
(I_{\bullet}\bullet L){\buildrel{Id_U\bullet (c_L\bullet
Id_L)}\over{\longrightarrow}} U\bullet((L o L^!)\bullet L)}
\end{eqnarray*}
We deduce that

\begin{eqnarray*}
(u')"=\matrix{U\bullet L\rightarrow (U\bullet I_{\bullet})\bullet L
{\buildrel{c_{\bullet}(U,I_{\bullet},L)}\over{\longrightarrow}}U\bullet
(I_{\bullet}\bullet L){\buildrel{Id_U\bullet (c_L\bullet
Id_L)}\over{\longrightarrow}} U\bullet((L o L^!)\bullet L)}
\end{eqnarray*}

\begin{eqnarray*}
\matrix{{\buildrel{Id_U\bullet{f_{L,L^!,L}}}\over{\longrightarrow}}U\bullet
(L o (L^!\bullet L)){\buildrel{Id_U \bullet (Id_L o
d_L)}\over{\longrightarrow}} U\bullet (L o I_o){\buildrel
{h_{U,L,I_o}}\over{\longrightarrow}}(U\bullet L) o
I_o{\buildrel{u}\over{\longrightarrow}} N}
\end{eqnarray*}
 Applying $(2.3)$ we deduce
 that

 \begin{eqnarray*}
 \matrix{U \bullet L\rightarrow U\bullet
(I_{\bullet}\bullet L){\buildrel{Id_U\bullet (c_L\bullet
Id_L)}\over{\longrightarrow}} U\bullet((L o L^!)\bullet L)}
\end{eqnarray*}

\begin{eqnarray*}
\matrix{
{\buildrel{Id_U\bullet{f_{L,L^!,L}}}\over{\longrightarrow}}U\bullet
(L o (L^!\bullet L)){\buildrel{Id_U \bullet (Id_L o
d_L)}\over{\longrightarrow}} U\bullet (L o I_o)\longrightarrow
U\bullet L}
\end{eqnarray*}

\begin{eqnarray*}
= Id_{U\bullet L}: U\bullet  L\longrightarrow U\bullet L
\end{eqnarray*}

This implies that:

\begin{eqnarray*}
(u')"=\matrix{U\bullet L{\buildrel {u}\over{\longrightarrow}} N}
\end{eqnarray*}

We have to show now that $(v")'=v$. We have:

\begin{eqnarray*}
(v")'=\matrix{U\longrightarrow U\bullet
I_{\bullet}{\buildrel{Id_U\bullet c_L}\over{\longrightarrow}}
U\bullet (L o L^!){\buildrel
{h_{U,L,L^!}}\over{\longrightarrow}}}(U\bullet L) o L^!\rightarrow
\end{eqnarray*}

\begin{eqnarray*}
((U\bullet L){\buildrel {v\bullet Id_L}\over{\longrightarrow}} (N
o L^!)\bullet L {\buildrel{f_{N,L^!,L}}\over{\longrightarrow}} N o
(L^!\bullet L){\buildrel{Id_N\bullet d_L}\over{\longrightarrow
N}})o L^!\longrightarrow N o L^!
\end{eqnarray*}

Using the fact that $o$ is a tensor product, we obtain that:

\begin{eqnarray*}
\matrix{((U\bullet L){\buildrel {v\bullet
Id_L}\over{\longrightarrow}} (N o L^!)\bullet L
{\buildrel{f_{N,L^!,L}}\over{\longrightarrow}} N o (L^!\bullet
L){\buildrel{Id_N\bullet d_L}\over{\longrightarrow N}})o L^!}
\end{eqnarray*}

\begin{eqnarray*}
=\matrix{(U\bullet L) o L^!{\buildrel {(v\bullet Id_L) o
Id_{L^!}}\over{\longrightarrow}} ((N o L^!)\bullet L) o L^!
{\buildrel{f_{N,L^!,L} o Id_{L^!}}\over{\longrightarrow}} (N o
(L^!\bullet L)) o L^!{\buildrel{(Id_N o d_L) o
Id_{L^!}}\over{\longrightarrow}} No L^!}
\end{eqnarray*}

This implies that

\begin{eqnarray*}
(v")'=\matrix{U\longrightarrow U\bullet
I_{\bullet}{\buildrel{Id_U\bullet c_L}\over{\longrightarrow}}
U\bullet (L o L^!){\buildrel {h_{U,L,L^!}}\over{\longrightarrow}}}
\end{eqnarray*}

\begin{eqnarray*}
\matrix{(U\bullet L) o L^!{\buildrel {(v\bullet Id_L) o
Id_{L^!}}\over{\longrightarrow}} ((N o L^!)\bullet L) o L^!
{\buildrel{f_{N,L^!,L} o Id_{L^!}}\over{\longrightarrow}} (N o
(L^!\bullet L)) o L^!{\buildrel{(Id_N o d_L) o Id_{L^!}}\over
{\longrightarrow}} (N o I_o) o
L^!{\buildrel{c_o(N,I_o,L^!)}\over{\longrightarrow}}N o (I_o o
L^!)\rightarrow N o L^!}
\end{eqnarray*}

Using property $(2.5)$ we obtain that:

\begin{eqnarray*}
\matrix{U\bullet (L o L^!){\buildrel
{h_{U,L,L^!}}\over{\longrightarrow}}(U\bullet L) o L^!{\buildrel
{(v\bullet Id_L) o Id_{L^!}}\over{\longrightarrow}} ((N o
L^!)\bullet L) o L^!}
\end{eqnarray*}

=
\begin{eqnarray*}
\matrix{U\bullet (L o L^!){\buildrel{v\bullet(Id_L o
Id_{L^!})}\over{\longrightarrow}}(N o L^!)\bullet (L o
L^!){\buildrel{h_{N o L^!,L,L^!}}\over{\longrightarrow}}((N o
L^!)\bullet L) o L^!}
\end{eqnarray*}

Thus

\begin{eqnarray*}
\matrix{U\longrightarrow U\bullet I_{\bullet}{\buildrel{Id_U\bullet
c_L}\over{\longrightarrow}} U\bullet (L o L^!){\buildrel
{h_{U,L,L^!}}\over{\longrightarrow}}(U\bullet L) o L^!{\buildrel
{(v\bullet Id_L) o Id_{L^!}}\over{\longrightarrow}} ((N o
L^!)\bullet L) o L^!}
\end{eqnarray*}

=

\begin{eqnarray*}
\matrix{U\longrightarrow U\bullet I_{\bullet}{\buildrel{Id_U\bullet
c_L}\over{\longrightarrow}}U\bullet (L o L^!){\buildrel{v\bullet
(Id_L o Id_{L^!})}\over{\longrightarrow}}(N o L^!)\bullet (L o
L^!){\buildrel{h_{N o L^!,L,L^!}}\over{\longrightarrow}}((N o
L^!)\bullet L) o L^!}
\end{eqnarray*}

Since $o$ is a tensor product, $Id_{L o L^!}=Id_L o Id_{L^!}$. Thus
using the fact that $\bullet$ is a tensor product, we have:

\begin{eqnarray*}
\matrix{U\bullet I_{\bullet}{\buildrel{Id_U\bullet
c_L}\over{\longrightarrow}}U\bullet (L o L^!){\buildrel{v\bullet
(Id_L o Id_{L^!})}\over {\longrightarrow}}(N o L^!)\bullet (L o
L^!)}
\end{eqnarray*}

=

\begin{eqnarray*}
\matrix{U\bullet I_{\bullet}{\buildrel {v\bullet
Id_{I_{\bullet}}}\over{\longrightarrow}}(N o L^!)\bullet
I_{\bullet}{\buildrel{Id_{N o L^!}\bullet
c_L}\over{\longrightarrow}}(N o L^!)\bullet (L o L^!)}
\end{eqnarray*}

This implies that

 \begin{eqnarray*} \matrix{U\longrightarrow
U\bullet I_{\bullet}{\buildrel{Id_U\bullet
c_L}\over{\longrightarrow}} U\bullet (L o L^!){\buildrel
{h_{U,L,L^!}}\over{\longrightarrow}}(U\bullet L) o L^!{\buildrel
{(v\bullet Id_L) o Id_{L^!}}\over{\longrightarrow}} ((N o
L^!)\bullet L) o L^!}
\end{eqnarray*}

=

\begin{eqnarray*}
\matrix{U{\buildrel v\over {\longrightarrow}} (N o L^!) \bullet
I_{\bullet}{\buildrel{Id_{N o L^!}\bullet
c_L}\over{\longrightarrow}} (N o L^!)\bullet (L o L^!){\buildrel
{h_{N o L^!,L,L^!}}\over{\longrightarrow}}((N o L^!)\bullet L) o
L^!}
\end{eqnarray*}

Using the fact that the associative constraint $c_o$ is a functorial
isomorphism, we deduce that

\begin{eqnarray*}
(v")'=\matrix{U{\buildrel v\over {\longrightarrow}} (N o L^!)
\bullet I_{\bullet}{\buildrel{Id_{N o L^!}\bullet
c_L}\over{\longrightarrow}} (N o L^!)\bullet (L o L^!){\buildrel
{h_{N o L^!,L,L^!}}\over{\longrightarrow}}((N o L^!)\bullet L) o
L^!}
\end{eqnarray*}

\begin{eqnarray*}
 \matrix{{\buildrel{f_{N,L^!,L} o
Id_{L^!}}\over{\longrightarrow}} (N o (L^!\bullet L)) o
L^!{\buildrel{{c_o}(N,L^!\bullet L,L^!)}\over{\longrightarrow}}N o
((L^!\bullet L)o L^!){\buildrel{Id_N o (d_L o
Id_{L^!})}\over{\longrightarrow}} No L^!}
\end{eqnarray*}

The property $(2.2)$ implies that:

\begin{eqnarray*}
\matrix{(N o L^!)\bullet (L o L^!){\buildrel{h_{N o
L^!,L,L^!}}\over{\longrightarrow}}((N o L^!)\bullet L) o L^!
{\buildrel{f_{N,L^!,L} o Id_{L^!}}\over{\longrightarrow}} (N o
(L^!\bullet L)) o L^! {\buildrel{{c_o}_{N,L^!\bullet
L,L^!}}\over{\longrightarrow}}N o ((L\bullet L^!) o L^!)}
\end{eqnarray*}

\begin{eqnarray*}
=\matrix{ (N o L^!)\bullet (L o L^!) {\buildrel{f_{N,L^!,L o
L^!}}\over{\longrightarrow}} N o (L^!\bullet (L o
L^!)){\buildrel{Id_N o h_{L^!,L,L^!}}\over{\longrightarrow}}N o
((L^!\bullet L) o L^!)}
\end{eqnarray*}

This implies that:
\begin{eqnarray*}
(v")'=\matrix{U{\buildrel v\over {\longrightarrow}} (N o L^!)
\bullet I_{\bullet}{\buildrel{Id_{N o L^!}\bullet
c_L}\over{\longrightarrow}} (N o L^!)\bullet (L o L^!)}
\end{eqnarray*}

\begin{eqnarray*}
\matrix{{\buildrel{f_{N,L^!,L o L^!}}\over{\longrightarrow}} N o
(L^!\bullet (L o L^!)){\buildrel{Id_N o
h_{L^!,L,L^!}}\over{\longrightarrow}}N o ((L^!\bullet L) o
L^!){\buildrel{Id_N  o(d_L o Id_{L^!})}\over{\longrightarrow}} N
oL^!}
\end{eqnarray*}

The property $(2.5)$ implies that:

\begin{eqnarray*}
\matrix{(N o L^!) \bullet I_{\bullet}{\buildrel{Id_{N o
L^!}\bullet c_L}\over{\longrightarrow}} (N o L^!)\bullet (L o L^!)
{\buildrel{f_{N,L^!,L o L^!}}\over{\longrightarrow}} N o
(L^!\bullet (L o L^!))}
\end{eqnarray*}

\begin{eqnarray*}
=\matrix{(N o L^!) \bullet
I_{\bullet}{\buildrel{f_{N,L,I_{\bullet}}}\over{\longrightarrow}}
N  o (L^!\bullet I_{\bullet}){\buildrel{Id_N o(Id_{L^!}\bullet
c_L)}\over{\longrightarrow}} N o (L ^!\bullet(L o L^!))}
\end{eqnarray*}

This implies that:

\begin{eqnarray*}
(v")'=\matrix{U{\buildrel v\over {\longrightarrow}}(N o L^!)
\bullet
I_{\bullet}{\buildrel{f_{N,L,I_{\bullet}}}\over{\longrightarrow}}
N  o (L^!\bullet I_{\bullet})}
\end{eqnarray*}

\begin{eqnarray*}
\matrix{{\buildrel{Id_N o(Id_{L^!}\bullet
c_L)}\over{\longrightarrow}} N o (L^! \bullet (L o
L^!)){\buildrel{Id_N o h_{L^!,L,L^!}}\over{\longrightarrow}} N o
((L^!\bullet L) o L^!)\rightarrow N o L^!}
\end{eqnarray*}

The property $(2.4)$ implies that:

\begin{eqnarray*}
\matrix{L^!\longrightarrow L^! \bullet
I_{\bullet}{\buildrel{Id_{L^!}\bullet
c_L}\over{\longrightarrow}}L^!\bullet (L o
L^!){\buildrel{h_{L,L^!,L}}\over{\longrightarrow}}(L^!\bullet L) o
L^!\longrightarrow L^!}
\end{eqnarray*}

is the identity of $L^!$. This implies that $(v")'=v$ since $o$ is a
tensor product.
\end{proof}

The previous theorem implies that the tensor category $(C,\bullet)$
is endowed with an internal $\underline{Hom}(U,V)= V o L^!$. The
general properties of internal $Hom$ ( see [4] Definition 1.6 page
109) implies the existence of an isomorphism:

\begin{eqnarray*}
d_{U_1,U_2}:\underline{Hom}(U_1,U_2)\bullet U_1\rightarrow U_2
\end{eqnarray*}

The map $d_{U_1,U_2}$ is the composition:

\begin{eqnarray*}
\matrix{(U_2 o U_1^!)\bullet
U_1{\buildrel{f_{U_2,U^!_1,U_1}}\over{\longrightarrow}} U_2 o
(U_1^!\bullet U_1){\buildrel{Id_{U_2} o
d_{U_1}}\over{\longrightarrow}} U_2 o I_o\longrightarrow U_2}.
\end{eqnarray*}

There exists a map:

\begin{eqnarray*}
l_{U_1,U_2,U_3}:\underline{Hom}(U_2,U_3)\bullet\underline{Hom}(U_1,U_2)
\longrightarrow \underline{Hom}(U_1,U_3)
\end{eqnarray*}

The map $l_{U_1,U_2,U_3}$ is the composition:

\begin{eqnarray*}
\matrix{(U_3 o U_2^!)\bullet (U_2 o
U_1^!){\buildrel{f_{U_3,U_2^!,U_2 o U_1^!}}\over{\longrightarrow}}
U_3 o (U_2^!\bullet (U_2 o U_1^!)){\buildrel{Id_{U_3} o
h_{U_2^!,U_2,U_1}}\over{\longrightarrow}} U_3 o ((U_2^!\bullet
U_2) o U_1^!){\buildrel{Id_{U_3} o (d_{U_2} o
Id_{U_1^!})}\over{\longrightarrow}} U_3 o U_1^!}.
\end{eqnarray*}

\begin{proposition}

Let $(C,\bullet, o)$ be a quadratic category. There exists a
canonical isomorphism $(U\bullet V)^!\longrightarrow V^! o U^!$.
\end{proposition}

\begin{proof}

The general properties of tensor categories imply the existence of
an isomorphism (see [4] 1.6.3 page 110):
\begin{eqnarray*}
\underline{Hom}(U_1\bullet
U_2,U_3)\longrightarrow\underline{Hom}(U_1,\underline{Hom}(U_2,U_3))
\end{eqnarray*}

we obtain an  isomorphism:

\begin{eqnarray*}
U_3 o (U_1\bullet U_2)^!\longrightarrow (U_3 o U_2^!) o U_1^!
\end{eqnarray*}

Suppose that $U_3=I_o$, the previous isomorphism  induces an
isomorphism:

\begin{eqnarray*}
(U_1\bullet U_2)^!\longrightarrow U_2^! o U_1^!
\end{eqnarray*}
\end{proof}

Let $h:U\rightarrow U'$ be a map. We define the morphism
$h^!:{U'}^!\rightarrow U^!$ as follows:

\begin{eqnarray*}
\matrix{{U'}^!\rightarrow {U'}^!\bullet
I_{\bullet}{\buildrel{Id_{{U'}^!}\bullet
c_U}\over{\longrightarrow}}{U'}^!\bullet (U o U^!){\buildrel
{h_{{U'}^!,U,U^!}}\over{\longrightarrow}}({U'}^!\bullet U) o
U^!{\buildrel{(Id_{{U'}^!}\bullet h)o
Id_{U^!}}\over{\longrightarrow}} ({U'}^!\bullet U') o
U^!{\buildrel{d_{U'} o I_{U^!}}\over{\longrightarrow}} U^!}
\end{eqnarray*}

\begin{proposition}

 Let $(C,o,\bullet)$ be a quadratic category, for every object $U$
of $C$, $U^!$ is unique up to an isomorphism.
\end{proposition}

\begin{proof}

The object $U^!$ represents the contravariant functor $V\rightarrow
Hom_C(V\bullet U,I_o)$. This implies that this object is unique up
to an isomorphism.
\end{proof}

We denote by $\underline{hom}(U,V)$ the object
${\underline{Hom}(U,V)}^!$, the paragraph before proposition $3$,
implies the existence of a   map $l_{U_1,U_2,U_3}$. The dual of
$l_{U_1,U_2,U_3}$ defines a map:

\begin{eqnarray*}
l^!_{U_1,U_2,U_3}:\underline{hom}(U_1,U_3)\longrightarrow
\underline{hom}(U_1,U_2) o \underline{hom}(U_2,U_3)
\end{eqnarray*}

The morphism $l_{U,U,U}$ endows $\underline{Hom}(U,U)$ with a
product
$l_{U,U,U}=l_U:\underline{Hom}(U,U)\bullet\underline{Hom}(U,U)\rightarrow
\underline{Hom}(U,U)$, and the morphism $l^!_{U,U,U}={l^!}_U:
\underline{hom}(U,U)\longrightarrow \underline{hom}(U,U)
o\underline{hom}(U,U)$ endows $\underline{hom}(U,U)$ with a
coproduct. There exists a canonical product on
$\underline{hom}(U,U)$ defined by

\begin{eqnarray*}
\matrix{(U^!\bullet U) o (U^!\bullet U){\buildrel{Id_{U^!\bullet
U} o d_U}\over{\longrightarrow}} U^!\bullet U}
\end{eqnarray*}

\begin{proposition}

 The object $(\underline{Hom}(U,U), l_U)$ is an algebra, and the
object $(\underline{hom}(U,U), {l^!}_U)$, is a coalgebra. This means
that the following diagrams commute:

\begin{eqnarray*}
\matrix{\underline{Hom}(U,U)\bullet
(\underline{Hom}(U,U)\bullet\underline{Hom}(U,U))
{\buildrel{Id_{\underline{Hom}(U,U)}\bullet l_U}\over
{\longrightarrow}}
\underline{Hom}(U,U)\bullet\underline{Hom}(U,U)\cr\cr \downarrow
c_{\bullet}(\underline{Hom}(U,U),\underline{Hom}(U,U),\underline{Hom}(U,U))\downarrow
l_U\cr\cr (\underline{Hom}(U,U)\bullet\underline{Hom}(U,U))\bullet
\underline{Hom}(U,U){\buildrel{l_U\bullet
Id_{\underline{Hom}(U,U)}}\over
{\longrightarrow}}\underline{Hom}(U,U)\bullet\underline{Hom}(U,U)
{\buildrel{l_U}\over{\longrightarrow}} \underline{Hom}(U,U)}
\end{eqnarray*}

\begin{eqnarray*}
\matrix{\underline{Hom}(U,U)\bullet
I_{\bullet}{\buildrel{Id_{\underline{Hom}(U,U)}\bullet
c_U}\over{\longrightarrow}}&\underline{Hom}(U,U)\bullet
\underline{Hom}(U,U)\cr\cr \downarrow &\downarrow l_U\cr\cr
\underline{Hom}(U,U){\buildrel{Id_{\underline{Hom}(U,U)}}\over{\longrightarrow}}
&\underline{Hom}(U,U)}
\end{eqnarray*}

These two diagrams endow $\underline{Hom}(U,U)$ with the structure
of an algebra. The next two diagrams endow $\underline{hom}(U,U)$
with the structure of a coalgebra.

\begin{eqnarray*}
\matrix{\underline{hom}(U,U){\buildrel{l^!_U}\over{\longrightarrow}}
&\underline{hom}(U,U) o \underline{hom}(U,U)\cr\cr \downarrow l^!_U
&\downarrow l^!_U o Id_{\underline{hom}(U,U)}\cr\cr
&(\underline{hom}(U,U) o \underline{hom}(U,U)) o
\underline{hom}(U,U) \cr &\downarrow c_o\cr\cr\underline{hom}(U,U) o
\underline{hom}(U,U)&{\buildrel{Id_{\underline{hom}(U,U)} o
l^!_U}\over{\longrightarrow}} \underline{hom}(U,U) o
(\underline{hom}(U,U) o \underline{hom}(U,U))}
\end{eqnarray*}

\begin{eqnarray*}
\matrix{\underline{hom}(U,U){\buildrel{l^!_U}\over{\longrightarrow}}
&\underline{hom}(U,U) o \underline{hom}(U,U)\cr\cr \downarrow
&\downarrow Id_{\underline{hom}(U,U)} o d_U\cr\cr
\underline{hom}(U,U)\longrightarrow &\underline{hom}(U,U)}
\end{eqnarray*}
\end{proposition}

\begin{proof}

We have only to show that the first two  diagrams are commutative
since the two last are their dual. The fact that first two diagrams
commute follows from general properties of tensor categories
\end{proof}

\begin{proposition}

 Let $(C,\bullet,o,!)$ be a quadratic category. The quadratic dual
$I_{\bullet}^!$ of $I_{\bullet}$ is isomorphic to $I_o$.
\end{proposition}

\begin{proof}

The theorem 2 implies that the functor defined on $C$ by
$U\rightarrow Hom_C(U,I_o)=Hom_C(U\bullet I_{\bullet},I_o)$ is
representable by $I_o\circ I_{\bullet}^!$ and $I_o$. Since $I_o\circ
I_{\bullet}^!$ is isomorphic to $I_{\bullet}^!$, we deduce that
$I_{\bullet}^!$ is isomorphic to $I_o$
\end{proof}

\begin{definition}

Let $(C, o,\bullet,!)$ be a quadratic category, the commutative
constraints $c'_{\bullet}$ of $(C,\bullet)$ and $c'_o$ of $(C,o)$
is a quadratic braiding if for every objects $U_1$, $U_2$ and
$U_3$ of $C$, the following diagram is commutative:

\begin{eqnarray*}
\matrix{U_1\bullet(U_2 o U_3)&{\buildrel{c'_{\bullet}(U_1,U_2 o
U_3)}\over{\longrightarrow}}(U_2 o U_3)\bullet
U_1{\buildrel{c'_o(U_2,U_3)\bullet
Id_{U_1}}\over{\longrightarrow}}&(U_3 o U_2)\bullet U_1\cr
\downarrow h_{U_1,U_2,U_3}&&\downarrow f_{U_3,U_2,U_1}\cr
(U_1\bullet U_2) o U_3 &{\buildrel{c'_o(U_1\bullet
U_2,U_3)}\over{\longrightarrow}}U_3 o (U_1\bullet
U_2){\buildrel{Id_{U_3} o
c'_{\bullet}(U_1,U_2)}\over{\longrightarrow}}& U_3 o(U_2\bullet
U_1)}
\end{eqnarray*}
\end{definition}

\begin{proposition}

 Let $(C,o,\bullet,!)$ be a quadratic category endowed with a
quadratic braiding, and $U$ an object of $C$, then ${U^!}^!$ is
isomorphic to $U$. In particular ${I^!_{\bullet}}^!\simeq
{I^!_o}\simeq I_{\bullet}$.
\end{proposition}

\begin{proof}

Let $U$ be an object of $C$, consider the morphisms

\begin{eqnarray*}
\matrix{c'_{U^!}:I_{\bullet}{\buildrel{c_U}\over{\longrightarrow}}U
o U^!{\buildrel{c'_o(U,U^!)}\over{\longrightarrow}}U^! o U}
\end{eqnarray*}

and

\begin{eqnarray*}
\matrix {d'_{U^!}:U\bullet
U^!{\buildrel{c'_{\bullet}(U,V)}\over{\longrightarrow}}U^!\bullet
U{\buildrel{d_U}\over{\longrightarrow}}I_o}
\end{eqnarray*}

We are going to show that the universal properties  verified by
$c_{U^!}$ and $d_{U^!}$ and dual of $U^!$ are verified by $c'_{U^!}$
and $d'_{U^!}$ and $U$.

We have:

\begin{eqnarray*}
\matrix{U^!\rightarrow I_{\bullet}\bullet
U^!{\buildrel{c'_{U^!}\bullet Id_{U^!}}\over{\longrightarrow}}(U^!
o U)\bullet U^!{\buildrel{f_{U^!,U,U^!}}\over{\longrightarrow}}
U^! o (U\bullet U^!){\buildrel{Id_{U^!} o
d'_{U^!}}\over{\longrightarrow}} U^! o I_o}
\end{eqnarray*}

\begin{eqnarray*}
=\matrix{U^!\longrightarrow I_{\bullet}\bullet
U^!{\buildrel{c_U\bullet Id_{U^!}}\over{\longrightarrow}} (U o
U^!)\bullet U^!{\buildrel{c'_o(U,U^!)\bullet
Id_{U^!}}\over{\longrightarrow}}(U^! o U)\bullet U^!}
\end{eqnarray*}

\begin{eqnarray*}
\matrix{    {\buildrel{f_{U^!,U,U^!}}\over{\longrightarrow}}U^! o
(U\bullet U^!){\buildrel{Id_{U^!} o
c'_{\bullet}(U,U^!)}\over{\longrightarrow}}U^! o (U^!\bullet
U){\buildrel{Id_{U^!} o d_U}\over{\longrightarrow}} U^! o I_o}
\end{eqnarray*}

The compatibility property of the quadratic braiding implies:

\begin{eqnarray*}
\matrix{ (U o U^!)\bullet U^!{\buildrel{c'_o(U,U^!)\bullet
Id_{U^!}}\over{\longrightarrow}}(U^! o U)\bullet U^!
{\buildrel{f_{U^!,U,U^!}}\over{\longrightarrow}}U^! o (U\bullet
U^!)}
\end{eqnarray*}

\begin{eqnarray*}
=\matrix{(U o U^!)\bullet U^!{\buildrel{{c'}_{\bullet}}(U o U^!,
U^!)\over{\longrightarrow}} U^!\bullet (U o
U^!){\buildrel{h_{U^!,U,U^!}}\over{\longrightarrow}}(U^!\bullet U) o
U^!{\buildrel{c'_o(U^!\bullet U,U^!)}\over{\longrightarrow}}U^! o
(U^!\bullet U){\buildrel{Id_{U^!} o
c'_{\bullet}(U^!,U)}\over{\rightarrow}} U^! o (U\bullet U^!)}
\end{eqnarray*}

Since the square of $c'_{\bullet}$ is the identity, this implies
that:
\begin{eqnarray*}
\matrix{U^!\rightarrow I_{\bullet}\bullet
U^!{\buildrel{c'_{U^!}\bullet Id_{U^!}}\over{\longrightarrow}}(U^!
o U)\bullet U^!{\buildrel{f_{U^!,U,U^!}}\over{\longrightarrow}}
U^! o (U\bullet U^!){\buildrel{Id_{U^!} o
d'_{U^!}}\over{\longrightarrow}} U^! o I_o}
\end{eqnarray*}

\begin{eqnarray*}
=\matrix{U^!\longrightarrow I_{\bullet}\bullet
U^!{\buildrel{c_U\bullet Id_{U^!}}\over{\longrightarrow}}(U o
U^!)\bullet U^!{\buildrel{{c'}_{\bullet}}(U o U^!,
U^!)\over{\longrightarrow}} U^!\bullet (U o
U^!){\buildrel{h_{U^!,U,U^!}}\over{\longrightarrow}}(U^!\bullet U)
o U^!{\buildrel{c'_o(U^!\bullet U,U^!)}\over{\longrightarrow}}U^!
o (U^!\bullet U)}
\end{eqnarray*}

\begin{eqnarray*}
\matrix{{\buildrel{Id_{U^!} o d_U}\over{\longrightarrow}} U^! o I_o}
\end{eqnarray*}

Since the commutative constraint of $\bullet$ is a functorial
isomorphism, we deduce that:

\begin{eqnarray*}
\matrix{U^!\longrightarrow I_{\bullet}\bullet
U^!{\buildrel{c_U\bullet Id_{U^!}}\over{\longrightarrow}}(U o
U^!)\bullet U^!{\buildrel{{c'}_{\bullet}(U o
U^!,U^!)}\over{\longrightarrow}}U^!\bullet (U o U^!)}
\end{eqnarray*}

\begin{eqnarray*}
=\matrix{U^!\longrightarrow I_{\bullet}\bullet
U^!{\buildrel{c'_{\bullet}(I_{\bullet},U^!)}\over{\longrightarrow}}U^!\bullet
I_{\bullet}{\buildrel{Id_{U^!}\bullet
c_U}\over{\longrightarrow}}U^!\bullet (U o U^!)=U^!\longrightarrow
U^!\bullet I_{\bullet}{\buildrel{Id_{U^!}\bullet
c_U}\over{\longrightarrow}} U^!\bullet (U o U^!) }
\end{eqnarray*}

This implies that

\begin{eqnarray*}
\matrix{U^!\rightarrow I_{\bullet}\bullet
U^!{\buildrel{c'_{U^!}\bullet Id_{U^!}}\over{\longrightarrow}}(U^!
o U)\bullet U^!{\buildrel{f_{U^!,U,U^!}}\over{\longrightarrow}}
U^! o (U\bullet U^!){\buildrel{Id_{U^!} o
d'_{U^!}}\over{\longrightarrow}} U^! o I_o}
\end{eqnarray*}

\begin{eqnarray*}
=\matrix{U^!\longrightarrow U^!\bullet
I_{\bullet}{\buildrel{Id_{U^!}\bullet
c_U}\over{\longrightarrow}}U^!\bullet (U o
U^!){\buildrel{h_{U^!,U,U^!}}\over{\longrightarrow}}(U^!\bullet U)
o U^!{\buildrel{c'_o(U^!\bullet U,U^!)}\over{\longrightarrow}}U^!
o (U^!\bullet U)}
\end{eqnarray*}

\begin{eqnarray*}
\matrix{\longrightarrow U^! o (U^!\bullet U){\buildrel{Id_{U^!} o
d_U}\over{\longrightarrow}} U^! o I_o}
\end{eqnarray*}

Using the fact that $c'_o$ is a quadratic braiding, we obtain that
$(U^!\bullet U) o U^!{\buildrel{c'_o(U^!\bullet
U,U^!)}\over{\longrightarrow}}U^! o (U^!\bullet
U){\buildrel{Id_{U^!} o d_{U}}\over{\longrightarrow}} U^! o
I_o=(U^!\bullet U) o U{\buildrel{d_U o
Id_{U^!}}\over{\longrightarrow}}I_o o
U^!{\buildrel{c'_o(I_o,U^!)}\over{\longrightarrow}}U^! o I_o$, since
$I_o o U^!{\buildrel{c'_o(I_o,U^!)}\over{\longrightarrow}}U^! o
I_o\rightarrow U^!=I_o o U^!\rightarrow U^!$  property 2.4, implies
 that

\begin{eqnarray*}
\matrix{U^!\rightarrow I_{\bullet}\bullet
U^!{\buildrel{c'_{U^!}\bullet Id_{U^!}}\over{\longrightarrow}}(U^!
o U)\bullet U^!{\buildrel{f_{U^!,U,U^!}}\over{\longrightarrow}}
U^! o (U\bullet U^!){\buildrel{Id_{U^!} o
d'_{U^!}}\over{\longrightarrow}} U^! o I_o}
\end{eqnarray*}

is the identity of $U^!$.

We also have:

\begin{eqnarray*}
\matrix{U\longrightarrow U\bullet I_{\bullet}{\buildrel{Id_U\bullet
c'_{U^!}}\over{\longrightarrow}}U\bullet(U^! o
U){\buildrel{h_{U,U^!,U}}\over{\longrightarrow}} (U\bullet U^!) o
U{\buildrel{d'_{U^!} o Id_U}\over{\longrightarrow}}I_o o U
\rightarrow U}
\end{eqnarray*}

\begin{eqnarray*}
=\matrix{U\longrightarrow U\bullet
I_{\bullet}{\buildrel{Id_U\bullet c_U}\over{\longrightarrow}}
U\bullet(U o U^!){\buildrel{Id_U\bullet c'_o(U,U^!)}\over
{\longrightarrow}} U \bullet (U^! o U)}
\end{eqnarray*}

\begin{eqnarray*}
\matrix{{\buildrel{h_{U,U^!,U}}\over{\longrightarrow}}(U\bullet U^!)
o U{\buildrel{c'_{\bullet}(U,U^!) o
Id_U}\over{\longrightarrow}}(U^!\bullet U) o U{\buildrel{d_U o
Id_U}\over{\longrightarrow}} I_o o U \rightarrow U}
\end{eqnarray*}

Using the compatibility property of the braiding, we have:

$h_{U,U^!,U}:U \bullet (U^! o U)\longrightarrow (U \bullet U^!) o
U $

\begin{eqnarray*}
=\matrix{U \bullet (U^! o U){\buildrel{c'_{\bullet}(U, U^! o
U)}\over{\longrightarrow}} (U^! o U)\bullet
U{\buildrel{c'_o(U^!,U)\bullet Id_U}\over{\longrightarrow}}(U o
U^!)\bullet U{\buildrel{f_{U,U^!,U}}\over{\longrightarrow}}U o
(U^!\bullet U){\buildrel{Id_U o
c'_{\bullet}(U^!,U)}\over{\longrightarrow}}U o(U\bullet
U^!){\buildrel{c'_o(U,U\bullet U^!)}\over{\rightarrow}} (U\bullet
U^!) o U}
\end{eqnarray*}
Using the fact that $c'_o$ is a functorial isomorphism and the
square of $c'_{\bullet}$ is the identity, this implies that:

\begin{eqnarray*}
\matrix{U\longrightarrow U\bullet I_{\bullet}{\buildrel{Id_U\bullet
c'_{U^!}}\over{\longrightarrow}}U\bullet(U^! o
U){\buildrel{h_{U,U^!,U}}\over{\longrightarrow}} (U\bullet U^!) o
U{\buildrel{d'_{U^!} o Id_U}\over{\longrightarrow}}I_o o U
\rightarrow U}
\end{eqnarray*}

\begin{eqnarray*}
=\matrix{U\longrightarrow U\bullet I_{\bullet}{\buildrel{Id_U\bullet
c_U}\over{\longrightarrow}} U\bullet(U o U^!){\buildrel{Id_U\bullet
c'_o(U,U^!)}\over {\longrightarrow}} U \bullet (U^! o U)}
\end{eqnarray*}

\begin{eqnarray*}
\matrix{{\buildrel{c'_{\bullet}(U, U^! o U)}\over{\longrightarrow}}
(U^! o U)\bullet U{\buildrel{c'_o(U^!,U)\bullet
Id_U}\over{\longrightarrow}}(U o U^!)\bullet
U{\buildrel{f_{U,U^!,U}}\over{\longrightarrow}}U o (U^!\bullet
U){\buildrel{c'_o(U,U^!\bullet U)}\over{\longrightarrow}}(U^!\bullet
U) o U{\buildrel{d_U o Id_U}\over{\longrightarrow}}I_o o U }
\end{eqnarray*}

Using the fact that $\bullet$ is a braided tensor product and the
fact that the square of $c'_o$ is the identity, we deduce that:

\begin{eqnarray*}
\matrix{U\longrightarrow U\bullet I_{\bullet}{\buildrel{Id_U\bullet
c_U}\over{\longrightarrow}} U\bullet(U o U^!){\buildrel{Id_U\bullet
c'_o(U,U^!)}\over{\longrightarrow}} U\bullet(U^! o U)
{\buildrel{c'_{\bullet}(U,U^! o U)}\over {\longrightarrow}} (U^! o
U)\bullet U{\buildrel{c'_o(U^!,U)\bullet
Id_U}\over{\longrightarrow}}(U o U^!)\bullet U}
\end{eqnarray*}

\begin{eqnarray*}
=\matrix{U\longrightarrow I_{\bullet}\bullet
U{\buildrel{c'_{\bullet}(I_{\bullet},U)}\over{\longrightarrow}}U\bullet
I_{\bullet}{\buildrel{Id_U\bullet c_U}\over{\longrightarrow}}
U\bullet (U o U^!){\buildrel{c'_{\bullet}(U,U o
U^!)}\over{\longrightarrow}}(U o U^!)\bullet U=U\rightarrow
I_{\bullet} \bullet U{\buildrel{c_U\bullet
Id_U}\over{\longrightarrow}}U\bullet (U o U^!)}
\end{eqnarray*}

This implies that:

\begin{eqnarray*}
\matrix{U\longrightarrow I_{\bullet} U\bullet
I_{\bullet}{\buildrel{Id_U\bullet
c'_{U^!}}\over{\longrightarrow}}U\bullet(U^! o
U){\buildrel{h_{U,U^!,U}}\over{\longrightarrow}} (U\bullet U^!) o
U{\buildrel{d'_{U^!} o Id_U}\over{\longrightarrow}}U o I_o}
\end{eqnarray*}

\begin{eqnarray*}
=\matrix{U\longrightarrow I_{\bullet}\bullet U{\buildrel{c_U\bullet
Id_U}\over{\longrightarrow}}(U o U^!)\bullet U
{\buildrel{f_{U,U^!,U}}\over{\longrightarrow}}U o (U^!\bullet
U){\buildrel{c'_o(U,U^!\bullet U)}\over{\longrightarrow}}(U^!\bullet
U) o U{\buildrel{d_U o Id_U}\over{\longrightarrow}} I_o o U }
\end{eqnarray*}

Using property $2.3$, and the fact that $o$ is a braided tensor
product we deduce that

\begin{eqnarray*}
\matrix{U\longrightarrow U\bullet
I_{\bullet}{\buildrel{Id_U\bullet
c'_{U^!}}\over{\longrightarrow}}U\bullet(U^! o
U){\buildrel{h_{U,U^!,U}}\over{\longrightarrow}} (U\bullet U^!) o
U{\buildrel{d'_{U^!} o Id_U}\over{\longrightarrow}}U o
I_o\rightarrow U}
\end{eqnarray*}

is the identity.

Since the quadratic dual is unique up to an isomorphism, we deduce
that $U$ is a quadratic dual of $U^!$.

\end{proof}

\begin{proposition}

 Let $(C,o,\bullet,!)$ be a quadratic category endowed with a
quadratic braiding. The contravariant functor $P:U\rightarrow U^!$
is an equivalence of category.
\end{proposition}

\begin{proof}

Let $U_1$, and $U_2$ be objects of $C$. Using the previous result,
we can suppose that ${U^!_1}^!=U_1$ and ${U_2^!}^!=U_2$.  The
theorem $2.2$, implies the existence of a bijection between
$Hom_C(U_1\bullet{U_2^!},I_o)$ and $Hom_C(U_1,U_2)$, and the
existence of an isomorphism between $Hom_C(U^!_2\bullet U_1, I_o)$
and $Hom_C(U^!_2,U^!_1)$. The commutative constraint $c'_{\bullet}$
defines an isomorphism between $Hom_C(U_1\bullet U^!_2,I_o)$ and
$Hom_C(U^!_2\bullet U_1,I_o)$. The object $U_1^!$ represents the
functor $n_{U_1}:V\rightarrow Hom_C(V\bullet U_1,I_o)$. This implies
that the correspondence defined on $C$ by $U\rightarrow U^!$ is
functorial, and using the Yoneda lemma, we deduce that the morphisms
between $n_{U_2}$ and $n_{U_1}$ are given by $Hom_C(U_2^!\bullet
U_1,I_o)$. This implies that $P$ is fully faithful.
\end{proof}

\begin{definition}

Let $(C,o,\bullet,!)$ be a quadratic category endowed with a
quadratic braiding. We denote by $C'$ the subcategory of $C$, such
that for each object $U$ of $C'$, there exists an object $V$ of $C$
such that $U$ is isomorphic to $V o I_{\bullet}$, $(C',o)$ is a
subtensor category of $(C,o)$. We say that the quadratic category
$(C,o,\bullet, !)$ is quadratic rigid if for every objects, $U_1$
and $U_2$ of $C'$, $U_1 o U_2=U_1\bullet U_2$ and the restriction of
the braided associativity constraints $f_{U_1,U_2,U_3}$ and
$h_{U_1,U_2,U_3}$ to $C'$ coincide with the associativity constraint
of $o$ and $\bullet$
\end{definition}

Let $(C,o,\bullet,!)$ be a rigid quadratic category, and  $U$ an
object of $C'$, we define $U^*$ to be $U^! o I_{\bullet}$. We have
$I_{\bullet}=I_o o I_{\bullet} $. This implies that $I_{\bullet}$ is
an object of $C'$. Thus $I_{\bullet}\bullet I_{\bullet}=I_{\bullet}
o I_{\bullet}=I_{\bullet}$. We deduce that for each object $U$ of
$C'$, $U o U^!$ is isomorphic to $(U o I_{\bullet}) o U^!$ which is
isomorphic to $U o (U^! o I_{\bullet})$. Since $U$ is  an object of
$C'$, $U o U^!$ is isomorphic to $U\bullet (U^! o
I_{\bullet})=U\bullet U^*$.

\begin{proposition}

 Let $U$ be an object of $C'$ the map $c'_U$, and the map

\begin{eqnarray*}
\matrix{d^1_{U^!}:U o(U^! o I_{\bullet}) =U\bullet (U^! o
I_{\bullet}){\buildrel{h_{U,U^!,I_{\bullet}}}\over{\longrightarrow}}(U
\bullet U^!) o I_{\bullet}{\buildrel{d'_U o
Id_{I_{\bullet}}}\over{\longrightarrow}} I_o o
I_{\bullet}=I_{\bullet}}
\end{eqnarray*}

defines on $(C',\bullet)$ the structure of a rigid tensor category.
\end{proposition}

\begin{proof}

We have to show  that for each element $U$ in $C'$:

\begin{eqnarray*}
\matrix{U\rightarrow U\bullet I_{\bullet}{\buildrel{Id_U\bullet
c'_U}\over{\longrightarrow}}U\bullet(U^! o U)=U\bullet (U^*\bullet
U){\buildrel{h_{U,U^!,U}}\over{\longrightarrow}}(U\bullet U^*) o U
{\buildrel{Id_U\bullet d^1_U }\over{\longrightarrow}}U}
\end{eqnarray*}
is the identity of $U$, and

\begin{eqnarray*}
\matrix{U^*\rightarrow I_{\bullet}\bullet U^*{\buildrel{c'_U\bullet
Id_{U^*}}\over{\longrightarrow}}(U^! o U)\bullet U^*=(U^*\bullet
U)\bullet U^*{\buildrel{f_{U^!,U,U^*}}\over{\longrightarrow}}U^*
o(U\bullet U^*){\buildrel{Id_{U^*} o
d^1_U}\over{\longrightarrow}}U^*}
\end{eqnarray*}
is the identity. These assertions follows from the fact that $U$ is
the dual of $U^!$. In the second assertion, we multiply $(2.4)$
applied to $U^!$ by $I_{\bullet}$. This implies that $U^*$ is a dual
of $U$ since the category $C$ is braided we deduce from Deligne
Milne (see [4]) that the category $C'$ is rigid
\end{proof}

\begin{definition}

Let $(C,o,\bullet, !)$ be a quadratic rigid category, and
$h:U\rightarrow U$ a morphism. We can define $h o Id_{I_{\bullet}}:U
o I_{\bullet}\rightarrow U o I_{\bullet}$. We define $Trace(h)$ to
be the trace of $h o Id_{I_{\bullet}}$. This is the endomorphism of
$I_{\bullet}$ defined by:

\begin{eqnarray*}
\matrix{I_{\bullet}{\buildrel{c'_U}\over{\longrightarrow}} U^! o
U{\buildrel{Id_U^! o h}\over{\longrightarrow}}U^! o
U{\buildrel{d^1_U}\over{\longrightarrow}}I_{\bullet}}
\end{eqnarray*}

We denote by $rank(U)$ the trace of $Id_U$. Let $h,h':U\rightarrow
U$, $Trace(hh')= Trace(h)Trace(h')$.
\end{definition}

The ring $Hom(I_{\bullet},I_{\bullet})$ is commutative since
$I_{\bullet}$ is the neutral element of $(C,\bullet)$.

\begin{proposition}

 Let $(C,o,\bullet,!)$ be a quadratic rigid category, and $C"$ be
the subcategory of $C$ such that for every object $U$ of $C"$, there
exists an object $V$ of $C$ such that $U=V\bullet I_o$. The category
$(C",o,I_o)$ is a rigid tensor category.
\end{proposition}

\begin{proof}

The restriction of the contravariant functor $P:C\rightarrow C$
defined on $C$ by $U\rightarrow U^!$ to $C'$ defines an isomorphism
between the tensor categories  $C'$ and $C"$.

\end{proof}

\begin{definition}

 Let $(C,o,\bullet,!)$ a quadratic category, and $h:U\rightarrow
V$ a morphism. The morphism $h':U^!\rightarrow V^!$ is a
contragredient of $h$ if and only if $h o h'\circ c_U=c_V$, and
$d_V\circ (h'\bullet h)=d_U$.
\end{definition}

\begin{proposition}

 Let $(C,o,\bullet,!)$ be a braided quadratic category. Then a
 map $h$ which has a contragredient is invertible.
\end{proposition}

\begin{proof}

Let $h':U^!\rightarrow V^!$ be a contragredient of the map
$h:U\rightarrow V$. We are going to show that $h'\circ h^!=Id_{V^!}$
and $h^!\circ h'=Id_{U^!}$, since the category is braided, this
implies that $h$ is invertible.

\begin{eqnarray*}
\matrix{h^!\circ
h'=U^!{\buildrel{h'}\over{\longrightarrow}}V^!{\buildrel{Id_{V^!}o
c_U}\over{\longrightarrow}}V^!\bullet(U o
U^!){\buildrel{h_{U^!,U,U^!}}\over{\longrightarrow}}(V^!\bullet U) o
U^!{\buildrel{(Id_{V^!}\bullet h) o
Id_{U}}\over{\longrightarrow}}(V^!\bullet V) o
U^!{\buildrel{d_V}\over{\longrightarrow}}U^!}
\end{eqnarray*}

Using property $2.5$, we deduce that $h^!\circ h'=$

\begin{eqnarray*}
\matrix{U^!{\buildrel{c_U}\over{\longrightarrow}}U^!\bullet(U o
U^!){\buildrel{h_{U^!,U,U^!}}\over{\longrightarrow}}(U^!\bullet U)
o U{\buildrel{h'\bullet h}\over{\longrightarrow}}(V^!\bullet V)  o
U^!{\buildrel{d_V}\over{\longrightarrow}} U^!}
\end{eqnarray*}

Using the fact that $d_V\circ (h'\bullet h)=d_U$, and property
$2.4$, we deduce that $h^!\circ h'=Id_{U^!}$. The proof that
$h^!\circ h'=Id_{U^!}$ is similar.

\end{proof}

\begin {definition}

A quadratic functor $H:(C,o,\bullet,!)\rightarrow
(C',o',\bullet',!')$ is a functor  of tensor categories
$H:(C,o)\rightarrow (C',o')$ which commutes with the quadratic dual,
 that is $H(U^!)=H(U)^!.$

\end{definition}

\begin{proposition}

 Let $(C,o,\bullet,!)$ and $(C',o',\bullet',!')$ be quadratic
braided categories, and $F,F':(C,o,\bullet,!)\longrightarrow
(C',o',\bullet',!')$ be two quadratic functors every morphism
$u:F\rightarrow F'$ is an isomorphism.
\end{proposition}

\begin{proof}

The morphism $u:F\rightarrow F'$ is defined by a family of morphisms
$u_U:F(U)\rightarrow F'(U)$ where $U$ is an object of $C$. The
morphism $u_{U^!}$ is a contragredient of $u_U$. This implies that
$u$ is an isomorphism. The inverse of $u$ is the family of maps
${u'_{U^!}}^!$
\end{proof}

We can generalize the morphisms defined by Manin [6] for quadratic
algebras as follows:

 We have the morphism $c_{I_o}:I_o\rightarrow
I_o o I^!_o=I_{\bullet}$.

Let $U$ and $V$ be objects of $C$, we have a morphism $U\bullet
V\longrightarrow U o V$ defined by

\begin{eqnarray*}
\matrix{U\bullet (V o I_o){\buildrel{c'_{\bullet}(U, V o
I_o)}\over{\longrightarrow}}(V o I_o)\bullet
U{\buildrel{f_{V,I_o,U}}\over{\longrightarrow}}V o (I_o\bullet
U){\buildrel{Id_V o(c_{I_o}\bullet Id_U)}\over{\longrightarrow}}V
o (I_{\bullet}\bullet
U){\buildrel{c'_o(V,U)}\over{\longrightarrow}}U o V}
\end{eqnarray*}

\section{Koszul complexes}

In this part, we suppose that the category $C$ is additive, and is
contained in an abelian category $C'$.
  Let $h$ be an element of $Hom(U,U)$, we denote by $d'_h$ the image
  of $h$ by the
isomorphism $Hom(U,U)\rightarrow Hom(I_{\bullet},U o U^!)$. We can
construct the map:

\begin{eqnarray*}
\matrix{d_h: U o U^!\longrightarrow I_{\bullet}\bullet (U o
U^!){\buildrel{ d'_h\bullet Id_{U o U^!}}\over{\longrightarrow}}(U o
U^!)\bullet (U o U^!){\buildrel{l_U}\over{\longrightarrow}}U o U^!}
\end{eqnarray*}

 \begin{definition}

 The category $C$ is an $n$-Koszul category, if for every object
 $U$, and each map $h:U\rightarrow U$, $(d_h)^{n+1}=0$. The
 category $C$ is a Koszul category if it is $1$-Koszul. This is equivalent
 to saying that  $(d_h)^2=0$.

\bigskip
Let $C$ be a Koszul category. We denote by $D_U$ the endomorphism
$d_{Id_U}$, and define the first Koszul complex to be $L(U)=(U o
U^!,D_U)$. We say that  $U$ is Koszul if $L(U)$ is exact. Let

\begin{eqnarray*}
\matrix{U{\buildrel{d_0}\over{\longrightarrow}}U_1...U_p
{\buildrel{d_p}\over{\longrightarrow}}U_{p+1}...}
\end{eqnarray*}

be a resolution of $U$ in $C'$. We say that this resolution is a
Koszul resolution, if there exists an object $V$ of $C$ endowed with
a differential $\alpha_V$, such that there exists  embedding
$\alpha_U: U o U^!\longrightarrow V$, $e_p:U_p\longrightarrow V$,
such that the following squares are commutative:

\begin{eqnarray*}
\matrix{U_p{\buildrel{d_p}\over{\longrightarrow}}U_{p+1}\cr
\downarrow e_p\downarrow e_{p+1}\cr
V{\buildrel{\alpha_V}\over{\longrightarrow}}V}
\end{eqnarray*}

\begin{eqnarray*}
\matrix{U o U^!{\buildrel{\alpha_U}\over{\longrightarrow}}V\cr
\downarrow D_U\downarrow \alpha_V\cr U o
U^!{\buildrel{\alpha_U}\over{\longrightarrow}}V}
\end{eqnarray*}
\end{definition}

\subsection{Koszul complexes of algebras}

We suppose that $C$ is a Koszul category, the objects of $C$ are
graded algebras defined over a field $F$. We denote by $U_i$ the
$i$-component of $U$, and we suppose that $U_0=F$,   the map $D_U$
is a left multiplication by an element $\alpha'_U$ of $U o U^!$, we
suppose that $C$ is stable by the usual tensor product of $F$-vector
spaces, and there exists an embedding $U o V\rightarrow U\otimes V$.
We denote by $\alpha_U$ the image of $\alpha'_U$ by this embedding.
We suppose that $\alpha_U\in U_1\otimes U^!_1$. The family
$(U\otimes U^!_i,\alpha_U)_{i\in{\N}}$ is a Koszul complex called
the first algebra Koszul complex. If this complex is exact, the
algebra is called a Koszul algebra.

\begin{theorem}

 Let $U$ be an object of $C$, if
 the complex $(U\otimes U^!,\alpha_U)$ is acyclic then
 $Ext_{U}(F,F)=U^!$.
\end{theorem}

\begin{proof}

 Suppose that $(U \otimes U^!,\alpha_U)$ is acyclic. This implies that:

\begin{eqnarray*}
0\longrightarrow U\longrightarrow U\otimes U^!_1...\longrightarrow
U\otimes {U^!}_i\longrightarrow U\otimes{U^!}_{i+1}...
\end{eqnarray*}

is a $U$-resolution of $U$. The tensor product $(U\otimes_F
{U^!}_i)\otimes_U F$ is isomorphic to $(U\otimes_U
F)\otimes_F{U^!}_i$. The tensor product $U\otimes_U F$ is $F$, since
the left-module of $F$ verifies $I(U)F=0$ where $I(U)$ is the
augmentation ideal. This implies that$(U\otimes_U
F)\otimes_F{U^!}_i={U^!}_i$ and the multiplication by $\alpha_U$
induces the zero map between $U^!_i$ and $U^!_{i+1}$. We deduce that
that $Ext_U(F,F)=U^!$$\bullet$

\end{proof}

\subsection{Second Koszul complex}

Let ${U^!}^*$ be the $F$-dual of $U^!$, $U\otimes {U^!}^*$ is
embedded in $Hom_U(U\otimes U^!,U)$ as follows: let $u_1,u_2$ be
elements of $U$, $v_1$ an element of ${U_l^!}^*$ and $v_2$ and
element of $U^!_p$. We define $(u_1\otimes v_1)(u_2\otimes
v_2)=v_1(v_2)u_1u_2$ if $l=p$, $(u_1\otimes v_1)(u_2\otimes v_2)=0$
if $l\neq p$. We can define the differential $D'_U$ on $U\otimes
{U^!}^*$ by setting $D'_U(h)(u)=h(\alpha_U(u))$. The complexes
$(U\otimes U^!,\alpha_U)$ and $(U\otimes {U^!}^*,D'_U)$ are dual
each other.

\subsection{ Quadratic algebras}

One of the main objective of specialists in Koszul structures is to
construct a Koszul resolution of an object $U$. The reason of this,
is the fact that using this resolution we  can  easily compute
$U$-homology and cohomology.
 The usual Bar complex allows to compute the homology of an algebra.
 Recall its definition. Let $U$ be an algebra, $\epsilon:U\rightarrow F$ the
augmentation of $U$, and $I(U)$ the kernel of $\epsilon$. The Bar
complex of $U$ is the tensor product $B(U)=U\otimes T(I(U))\otimes
U$. Its elements are denoted by $u[u_1:...:u_p]u'$. The Bar complex
is bigraded the degree of $v=u[u_1:...:u_p]u'$ is $(m,n)$ where
$m=\sum_{i=1}^{i=p}degree(u_i)+degree(u)+degree(u')$, and $n=p$ is
the homological degree.
 The converse of  theorem 3.2  for quadratic algebras is shown by Priddy [7] using a
 spectral sequence. Here is an elementary proof:

 \begin{proposition}

 Let $U$ be a quadratic algebra, if $Ext_U(F,F)=U^!$, then
 $U$ is a Koszul algebra.
\end{proposition}

\begin{proof}

Suppose that $U=T(V)/C$. We compute $Ext_U(F,F)$ using the bar
complex. We have $F\otimes B(U)\otimes F=T(I(U))$. This implies that
$Ext_U(F,F)=H^*(T(I(U))^*)$. We denote by $D'$ the differential of
this complex. The algebra $T(V^*)$ is contained in $T(I(U)^*)$ since
$I(U)$ contains $V$. We use the homological degree to define the
graduation of $T(I(U)^*)$. We obtain that $T(I(U)^*)_0=T(V^*)$, and
$T(I(U))_1=\sum (V^{\otimes^p}\otimes (T(V)/C)\otimes
V^{\otimes^l})^*$. The kernel of the restriction of $D'$ to $T(V^*)$
is $U^!$. The complex $(I(U)\otimes {U^!}^*,D'_U)$ is contained in
the Bar complex $T((I(U),D)$. The result follows from the fact that
$H^{p,l}(T(I(U))=H_{p,l}(T(I(U))=0$ if $l\neq p$.
\end{proof}

\section{Representations of quadratic categories}

In this section we study representations of quadratic categories to
the category of quadratic algebras, and adapt the results of [4].
\begin{definition}

Let $U=L[V]/C$ be a quadratic algebra. We denote by $Aut(U)$. The
automorphisms group of $U$. An element $h$ of $Aut(U)$ is defined
by an automorphism $h_1$ of $V$ such that $(h_1\otimes h_1)(C)=C$.
Consider an algebraic group $H$, a quadratic representation of
$H$, in $U$ is a morphism    $\rho:H\longrightarrow Aut(U)$.
\end{definition}

\begin{proposition}

 Let  $H$ an algebraic group, the
set of quadratic representations of $H$  is a quadratic category.
\end{proposition}

\begin{proof}

Let $\rho_1:H\longrightarrow Aut(U_1)$, and $\rho_2:H\longrightarrow
Aut(U_2)$ be two quadratic representations, $\rho_1 o\rho_2$ is the
representation defined as follows: for each $h\in H$  $(\rho_1 o
\rho_2)(h):U_1 o U_2\longrightarrow U_1 o U_2$ is the morphism
$\rho_1(h) o \rho_2(h)$.

The map $(\rho_1\bullet\rho_2)(h)$ is the automorphism of
$\rho_1(h)\bullet\rho_2(h)$ of $U_1\bullet U_2$.

The map $(\rho_1^!)(h)$ is the automorphism $(\rho_1(h))^!$ of
$U_1^!$.

The category $C_H$ of quadratic representations endowed with the
tensor products $o$ and $\bullet$ that we have just defined is a
quadratic category.

The respective neutral elements $I^H_o$ and $I^H_{\bullet}$ for $o$
and $\bullet$ are the respective  trivial representations
$H\rightarrow Aut(I_o)$ and $H\rightarrow Aut(I_{\bullet})$
\end{proof}

Let $V$ be a vector space, $V$ is a quadratic algebra endowed with
the zero product. In this case $C=V^{\otimes^2}$. Thus the
category of representations of $H$ is embedded in the category of
quadratic representations of $H$.

\begin{definition}

A morphism of the category $C_H$ of quadratic representations is
defined by a family of automorphisms $e_U$ of $U$, where
$\rho:H\longrightarrow Aut(U)$ is an object of $C_H$, $U=T(U_1)/C$
such that $e_U$ commutes  with the action of $H$ on $U_1$. We
suppose that:

$e_{U\bullet V}=e_U\bullet e_V$

$e_{U o V}=e_U o e_V$

$e_{U^!}={e_U}^!$.
\end{definition}

We denote $Aut(C_H)$ the group of automorphisms of $C_H$.

\begin{proposition}
 The natural embedding $H\rightarrow Aut(C_H)$ is an isomorphism.
 \end{proposition}

\begin{proof}

The proof follows from the corresponding result for tensor
categories since the category of quadratic representations contains
the category of representations (see [4]).

\end{proof}

\begin{theorem}

 Let $(C,o,\bullet,!)$ be a quadratic braided  abelian category
such that $End(I_{\bullet})$ is the ground field $L$, and $P$ the
category of quadratic algebras defined over the field $L$. Suppose
that there exists an exact faithful functor $F:C\longrightarrow P$,
 then $(C,o,\bullet,!)$ is equivalent to the category of quadratic
 representations of an affine group scheme.
\end{theorem}

\begin{proof}

Let $P'$ be the category whose objects are couples $(N,T)$ where $N$
is a finite dimensional vector space, and $T$ a subspace of
$N\otimes N$. A morphism $h:(N,T)\longrightarrow (N',T')$ is a
linear map $h:N\longrightarrow N'$ such that $(h\otimes h)(T)\subset
T'$ The category $P'$ is endowed with the tensor product defined by
$(N,T)\otimes (N',T')=(N\otimes N',t_{23}(T\otimes T'))$.  The
functor $D:P\rightarrow P'$, $U=T(U_1)/C\longrightarrow (U_1,C)$ is
fully faithful, and it is a tensor functor when $P$ is endowed with
the tensor product $\bullet_P$ defined by $T(U_1)/C\bullet_P
T(U'_1)/C'=T(U_1\otimes U'_1)/t_{23}(C\otimes C')$.

Let $(C,o,\bullet,!)$ be a quadratic braided abelian category
endowed with an exact faithful quadratic functor $F:C\rightarrow P$,
then $D\circ F:C\rightarrow P'$ is exact and faithful. In Deligne
Milne [4], it is shown that $C$ is equivalent to the category of
comodules over a coalgebra ${\cal H}$. Since the object of $P'$ are
finite dimensional spaces as in Deligne Milne [4], we can endow
${\cal H}$ with a structure of a bialgebra since every morphism
between two quadratic functors is an isomorphism, and deduce that
$(C,\bullet,o,!)$ is equivalent to the category or quadratic
representations of the affine group scheme $H$ defined by ${\cal H}$
\end{proof}

 In the proof of theorem 5, we can define on $P'$ the following
 tensor structure: Let $(U,T)$, and $(U',T')$ be two elements of
 $P'$, $(U,T) o (U',T')=(U\otimes U',t_{23}(U\otimes T'+ T\otimes
 U'))$. This defines on the coalgebra ${\cal H}$ another algebra
 product.

\begin{definition}

A neutral braided Tannakian quadratic category is a quadratic
category $(C,o,\bullet,!)$ such that there exists an exact faithful
$L$-linear functor $F:C\longrightarrow P$, where $P$ is the category
of quadratic algebras.   The functor $F$ is called a fibre functor.
If $L'$ is a $L$-algebra, a $L'$-fiber functor is a fiber functor
$F$ such that $D o F$ (see proof of theorem 5 for the definition of
D) takes its values in the category of $L'$-modules.

\end{definition}

\begin{theorem}

Let $(C,o,\bullet,!)$ be a Tannakian quadratic category, $L$ a field
and $Aff_L$ the category of affine  schemes defined over $L$. We
denote by $Fib(C)$ the category of fiber functors of $C$. The
functor $Fib(C)\rightarrow Aff_L$, $F\rightarrow H$ is a gerbe over
$Aff_L$ The fiber of $Spec(L')$ where $L'$ is a $L$-algebra is the
$L'$-valued fiber functor.
\end{theorem}

\begin{proof}

Let $F':C(o,\bullet,!)\rightarrow P$ be another fiber $L'$-functor,
we have to show that $Hom(F',F')$ is representable by a $H$-torsor
over $Spec(L')$ where $H$ represents $F$. The composition of a fiber
functor and the functor $D$ defined in the proof of theorem $5$
defines fiber a functor $C$ to the category of $L'$-module. We apply
the corresponding result for Tanakian category to $(F_0(C),\bullet)$
[4]
\end{proof}

\section{Quadratic motives}

 Let $L$ be a field, and $V_L$ the
category of projective schemes defined on $L$. A cohomological
functor $F:V\rightarrow (C,\otimes)$, where $(C,\otimes)$ is a
tensor abelian category,  which verifies standards properties that
verified cohomologies theories like Kunneth formula $F(U\times
U')=F(U)\otimes F(U')$. The category of motives $F_0:V_L\rightarrow
N_L$ is an initial object in the category of cohomological functors,
that is for every cohomological functor $F:V\rightarrow C$, there
exists a functor of abelian categories $F_C:N_L\rightarrow C$ such
that $F=F_C\circ F_0$. The construction of a category of motives is
given in  Deligne-Milne [4].

Let $U$ be a quadratic algebra. We denote by $Pro(U)$ the non
commutative projective scheme defined by $U$. This is the category
of $U$-graded modules up to the modules of finite length. The
category of coherent sheaves over every projective scheme is
equivalent to the non commutative projective scheme defined by a
quadratic algebra. This realization is not unique, that is the
category of coherent sheaves defined on a projective scheme can be
equivalent to the non commutative projective schemes defined by
two non isomorphic quadratic algebras.

We define the category of quadratic motives to be the category of
$L$-quadratic algebras this is an abelian category.

\end{document}